\magnification=\magstep1
\input amstex
\documentstyle{amsppt}

\NoBlackBoxes

\pagewidth{5.15 truein}
\pageheight{8.00 truein}

\define\F{\Cal F}
\define\emp{\varnothing}

\topmatter
\title
Noncommutative Algebras Associated to Complexes and Graphs
\endtitle

\author
Israel Gelfand, Sergei Gelfand, Vladimir Retakh
\endauthor

\address
Department of Mathematics, Rutgers University, 110 Frelinghysen Road,
Piscataway, NJ 08854-8019
\endaddress
       
\email
igelfand\@math.rutgers.edu
\endemail

\address
American Mathematical Society, P.O.Box 6248, Providence, RI 02940
{\rm and} Institute for Problems of Information Transmission, 19, Ermolova str.,
Moscow, 103051, Russia
\endaddress

\email
sxg\@ams.org
\endemail

\address
Department of Mathematics, Rutgers University, 110 Frelinghysen Road,
Piscataway, NJ 08854-8019
\endaddress

\email
vretakh\@math.rutgers.edu
\endemail

\endtopmatter

\document

\subhead 1. Introduction \endsubhead
This is a first of our papers devoted to "noncommutative topology and
graph theory". Its origin is the paper \cite {GRW} where a new class
of noncommutative algebras $Q_n$ was introduced.
As explained in \cite {GRW}, the algebra $Q_n$ is closely related
to decompositions of a generic polynomial $P(t)$ of degree $n$
over a division algebra into linear factors.

The structure of the algebra $Q_n$ seems to be very interesting. It
has linearly independent generators $u(B)$,
$B\subset \{1,\dots ,n\}$. Here $u(\emptyset )=1$ is a unit element of
$Q_n$. An important property of $Q_n$ is that
under any homomorphism of $Q_n$ into a commutative integral domain,
each element $u(B)$ with $|B|\geq 2$ maps to zero. In other words,
elements 
$u(B)$ carry the "noncommutative nature" of $Q_n$. Moreover, the
"degree of noncommutativity" carried by $u(B)$ depends on the size
of $B$.

The noncommutative nature of $Q_n$ can be studied by looking at
quotients of $Q_n$ by ideals generated by some $u(B)$. These 
quotients are "more commutative" then $Q_n$. For example, the
quotient of $Q_n$ by the ideal generated by all $u(B)$ with
$|B|\geq 2$ is
isomorphic to the algebra of commutative polynomials in $n$
variables.
 
To consider more refined cases we need to turn to a "noncommutative
combinatorial topology". In our approach the algebra $Q_n$
corresponds to an $n$-simplex $\Delta _n$ and we consider
quotients of $Q_n$ by ideals generated by some $u(B)$ 
corresponding to subcomplexes of $\Delta _n$. We describe generators and
relations for those quotients.

We pay special attention to the quotients of $Q_n$
corresponding to 1-dimensional subcomplexes of $\Delta _n$.
(They are "next" to algebras of commutative polynomials).

The third author was partially supported by the National
Science Foundation.

\subhead 2. The algebra $Q_n$ \endsubhead

For a natural $n$ by $I_n$ we denote the set $I_n= \{1,2,\dots,n\}$. A
noncommutative unital  algebra $Q_n$ (see \cite{GRW}) 
is the algebra generated by elements
$z_{A,i}$, where $A\subset I_n$, $i\notin A$, subject  to two groups of relations: 
additive relations 
$$
z_{A\cup i,j} + z_{A,i} = z_{A\cup j,i} + z_{A,j}
\tag 1
$$
and multiplicative relations
$$
z_{A\cup i,j} z_{A,i} = z_{A\cup j,i} z_{A,j}.
\tag 2
$$

In \cite{GRW} the authors
defined elements  $z_{A,B}\in Q_n$ for each pair of disjoint subsets
$A,B\subset I_n$. Of particular importance to us will be 
the elements $z_{\emp,A}$, $A\subset I_n$, for which we will use the notation
$u(A)$. These elements are defined as follows. Let $A\subset I_n$. Choose
$i\in A$. Then
$$
u(A) = \sum_{D\subset A,\,D\not\ni i} (-1)^{|A|-|D|} z_{D,i}.
\tag 3
$$
From additive relations (1) one can easily see that the right-hand side of (3)
does not depend on the choice  of $i\in A$.

We have the following simple result.

\proclaim{Lemma} Let $L\subset Q_n$ be the linear subspace spanned by all
$z_{A,i}$. Then the elements  $u(A)$, $A\subset I_n$, form a basis of $L$.
\endproclaim

\demo{Proof} Formula (3) implies that 
$$
z_{A,i} = \sum_{D\subset A} u(D\cup i).
$$
Therefore, the elements $u(A)$ generate $L$. On  the other hand, it is proved
in \cite{GRW} that $\dim L =  2^n -1$.\qed
\enddemo

\proclaim{Corollary} The algebra $Q_n$ can be defined as the algebra with
generators $u(A)$, $A\subset  I_n$, and relations
$$
\multline
\sum_{C,D\subset A}\left(u(C\cup j)+u(C\cup i\cup j)\right)u(D\cup i)\\
=\sum_{C,D\subset A}\left(u(D\cup i)+u(D\cup i\cup j)\right)u(C\cup j),
\endmultline
\tag{$4_{A,i,j}$}
$$
one for each triple $(A,i,j)$, $A\subset I_n$, $i,j\notin A$, $i \neq j$.
\endproclaim

We will use formulas $(4_{ A,i,j})$ in the following form:
$$
\multline
\sum_{C,D\subset A}\left[u(C\cup i),u(D\cup j)\right]
\\
=\Big(\sum_{E\subset A}u(E\cup i\cup j)\Big)\sum_{F\subset A}\left( u(F\cup i)
- u(F\cup j)\right). 
\endmultline
\tag{$5_ {A,i,j}$}
$$

\subhead 3. Complexes \endsubhead 
\definition{Definition} (i) A {\it complex\/} with $n$ nodes is a family $\F$
of nonempty subsets $A\subset I_n$ satisfying the {\it filtering\/} condition:
$$
A\in \F, B\subset A \implies B\in \F.
\tag 6
$$

(ii)  The {\it dimension \/} of a complex $\F$ is defined as follows:
$$
\dim \F = \max_{A\in \F} (|A| - 1).
$$
\enddefinition

With any complex $\F$ we can associate a triangulated topological space
$T(\F)$, called the {\it geometrical realization of $\F$,} as follows.

Let $\Delta_n$ be the standard $n$-dimensional simplex with vertices numbered
$1,\dots, n$. To each  $F\subset I_n$ there corresponds the
$(|F|-1)$-dimensional face $\Delta(F)$ of $\Delta_n$ (the vertex $i$ if 
$F=\{i\}$, the open $(d-1)$-dimensional simplex with vertices
$i_1,i_2,\dots,i_d$ if $F=\{ i_1,i_2, \dots,  i_d\}$, $d>1$). Define
$$
T(\F) = \bigcup_{F \in\F}\Delta(F)\subset \Delta_n.
$$
Condition (6) implies that $T(\F)$ is a closed subset of $\Delta_n$.

\subhead 4. The algebra $Q(\F)$ \endsubhead 
\definition{Definition} Let $\F$ be a complex with $n$ nodes. Define $Q_n(\F)$
as the quotient algebra of  $Q_n$ by the ideal generated by the elements $u(A)$
for all $A\notin \F$.
\enddefinition

\remark{Remarks} 1. If $\F$ is the family of all subsets of $I_n$, then
$Q_n(\F)= Q_n$. 

2. If $\dim \F=0$ then the algebra $Q(\F)$ is isomorphic to
the algebra of commutative polynomials in $n$ variables.

3. Let $n_1<n_2$ and let $\F$ be a complex with $n_1$ nodes. The inclusion
$I_{n_1}\subset I_{n_2}$  allows us to view $\F$ as a complex $\F'$ with
$n_2$ nodes. On the other hand, $Q_{n_1}$ is naturally 
isomorphic to the quotient of $Q_{n_2}$ by the ideal generated by all $z_{A,i}$
with $A\cup i\not\subset  I_{n_1}$. Formula (1) and the definition of $Q_n(\F)$
show that $Q_{n_1}(\F)$ is naturally isomorphic to 
$Q_{n_2}(\F')$. Therefore we can (and will) always assume that our complex $\F$
contains all one-elements subsets and will write $Q(\F)$ instead of $Q_n(\F)$. 

4. If $\F'\subset \F$ is a subcomplex, then $Q(\F')$  is naturally
isomorphic to
a quotient algebra of $Q(\F)$.
\endremark

The main result of this note is a description of the algebra $Q(\F)$ and its
properties in the case where $\F$  is a graph (i.e., a one-dimensional
complex). Before formulating the main theorem, we prove some  properties of the
algebras $Q(\F)$ for an arbitrary complex $\F$.

\proclaim{Proposition} Let $\F$ be a complex, $A,B\in \F$. Let there exist
$i\in A$, $j\in B$ such that 
$$
(i\cup A)\notin \F,\quad (j\cup B)\notin B.
\tag 7
$$ 
Then $[u(A),u(B)]=0$ in $Q(\F)$.
\endproclaim

\demo{Proof} We prove that for each pair $A'\subset A\setminus i$, $B'\subset
B\setminus j$ we have 
$$
[u(A'\cup i),u(B'\cup j)]=0
\tag 8
$$ 
in $Q(\F)$. For  $A'=A\setminus i$, $B'=B\setminus j$ we get the desired
result. The proof is by induction  in $|A'|+|B'|$. 

Before proceeding, let us note that condition (7) allows us to rewrite the
relation $(5_{A'\cup B',i,j)}$ in the  form
$$
\sum_{C\subset A',\,D\subset  B'}[u(C\cup i,D\cup j] = 0
\tag{$9_{A',B',i,j}$}
$$
Now we prove (8) by induction in  $|A'| + |B'|$.

If  $|A'| + |B'|=0$, i.e., $|A'| = |B'|=\emp$, the left-hand side of
$(9_{A',B',i,j})$ is reduced to a single term  
$[u(i),u(j)]$, and we get (8).

For an arbitrary $A'$, $B'$ each term in $(9_{A',B',i,j})$ except the term
$[u(A'\cup i,B'\cup j]$ is of the  form $[u(C\cup i,D\cup j]$ with $|C|+|D| <
|A'|+|B'|$. By the induction assumption, all these terms vanish 
and $(9_{A',B',i,j,})$ becomes $[u(A'\cup i),u(B'\cup j)]=0$. \qed
\enddemo

\subhead 5. The main theorem \endsubhead Now let $\F$ be a graph,
i.e., a complex of  dimension~1. Denote by $E=E(\F)$ the set of edges of $\F$,
i.e., the set of unordered pairs $(ij)$, $i\neq j$,  such that $\{i,j\}\in \F$.

\proclaim{Theorem} The algebra $A(\F)$ is generated by the elements $u(i)$,
$i\in I_n$, and $u(ij)$, $(ij)\in  E(\F)$. These elements satisfy the following relations 
(we assume that $u(ij)=0$ if $(ij)\notin E(\F)$.

(i)  For each pair $(i,j)$ of distinct integers $i,j$ between $1$ and $n$ we have 
$$
[u(i),u(j)]=u(ij)(u(i)-u(j)).
$$ 
In particular, $[u(i),u(j)]=0$ for $(ij)\notin E(\F)$.

(ii) For each triple $(i,j,k)$ of distinct integers 
between $1$ and $n$ we have 
$$
[u(ik),u(jk)] +[u(ik),u(j)] +[u(i),u(jk)] =u(ij)(u(ik) - u(jk)) .
$$.

(iii) For each quadruple $(i,j,k,l)$ of distinct integers
between $1$ and $n$ we have $[u(ij),u(kl)]=0$. 
\endproclaim

\demo{Proof}
First, we modify the relation $(5_{A,i,j})$ for the case
where $\F$ is a graph. Since in  this case $u(A)=0$ whenever $|A| \geq 3$, we
have, for each $A\subset I_n$,  
$$
\aligned
\sum_{k,l\in A}[u(ik) ,u(jl)] &+\sum_{k\in A} [u(ik),u(j)] + \sum_{k\in A} [u(i),u(jk)] + [u(i),u(j)]
\\ &- \sum_{k\in A}u(ij)(u(ik) - u(jk)) + u(ij)(u(i)-u(j))  = 0.
\endaligned
\tag{$10_{A,i,j}$}
$$
Denote the left-hand side of $(10_{A,i,j})$ by $R(A,i,j)$. If $A$ is empty, then 
$$
R(\empty,i,j)= [u(i),u(j)]-u(ij)(u(i)-u(j)).
$$
and formula ($10_{\empty,i,j}$) gives part (i) of the main theorem. Now we assume that 
$A$ is nonempty.  Choose $k\in A$. Then
$$
\multline
R(A,i,j) = R(A\setminus k,i,j) + [u(ik) ,u(jk)] 
+\sum_{l\in A\setminus k}  [(u(il), u(jk)] +\sum_{l\in A\setminus k} [(u(ik), u(jl)] 
\\+ [u(ik),u(j)] + [u(i), u(jk)]- u(ij)(u(ik) - u(jk)) 
\endmultline
\tag{11}
$$

Let is prove (ii). Let $i,j,k$ be three distinct elements of $I_n$. 
Take $A=\{k\}$. Since $R(\empty,i,j)=0$, relations $(10_{A,i,j})$ and (11) give 
$$
[u(ik) ,u(jk)] +  [u(ik),u(j)] + [u(i), u(jk)]- u(ij)(u(ik)- u(jk)) = 0,
\tag {12}
$$
i.e., (ii) is proved.

(iii)  Let $i,j,k,l$ be four distinct elements of $I_n$. Take $A=\{k,l\}$.
Since $R(\{l\},i,j)\allowmathbreak
=0$, we have from $(10_{A,i,j})$ and (11) that
$$
\multline
[u(ik) ,u(jk)]  + [u(ik),u(j)] + [u(i), u(jk)]-   u(ij)(u(ik) - u(jk))\\
 +   [ (u(ik), u(jl)]  - [ (u(jk), u(il)] = 0.
\endmultline
$$
The sum of the first four terms in the left-hand side of this equality vanishes
because of (12). Therefore, for  each ordered quadruple $(i,j,k,l)$ we have
$$
[ (u(ik), u(jl)]  = [ (u(jk), u(il)].
$$
Taking now the same four elements in the different order $(k,l,i,j)$ we obtain
the relation 
$$
[ (u(ik), u(jl)]  = [ (u(il), u(jk)].
$$
The last two relations imply that $[ (u(ik), u(jl)]  = 0$, i.e. (iii) is proved.

(iv) It remains to prove that relation $(10_{A,i,j})$ with $|A|\geq 3$ follows
from (i)--(iii). But, by (ii) and  (iii), for $k\in A$ relation (11) takes the
form $ R(A,i,j) = R(A\setminus k,i,j)$, and it remains to use induction on $|A|$.
\qed
\enddemo

\enddocument